\documentclass[10pt]{article}
\usepackage{latexsym,amsthm,amscd}
\newcounter{alphthm}
\newtheorem{thm}{Theorem}

\newtheorem{defn}{Definition}
\newtheorem{lem}{lemma}
\usepackage{amsmath}
\input{amssym}
\def\bd{\begin{defn}}
\def\ed{\end{defn}}
\def\bl{\begin{lem}}
\def\el{\end{lem}}
\def\bt{\begin{thm}}
\def\et{\end{thm}}

\newcommand{\R}{I\!\! R}
\numberwithin{equation}{section} \textwidth 14cm 

\begin{document}
\title{ Invariant conformal geometry on Finsler manifolds\thanks{The Fifth Conference of Balkan
 Society of Geometers , August 29 - September 2, 2005, Mangalia-
Romania, pp. 34-43.  \newline $\copyright$ Balkan Society of
Geometers, Geometry Balkan Press 2006.}}
\date{}
\author{B. Bidabad $\,\,\,\,\,\,\,\,\,\,\,$S. Hedayatian}
 \maketitle
\begin{abstract}
The electric capacity of a conductor in the $3$-dimensional
Euclidean space $I\!\!R^3 $ is defined as a ratio of a given
positive charge on the conductor to the value of potential on the
surface. This definition of the capacity is independent of the given
charge. The capacity of a set as a mathematical notion was defined
first by N. Wiener (1924) and was developed by O. Forstman
[\ref{R:FO}], C. J. de La Vallee Poussin, and several other French
mathematicians in connection with potential theory.
 This paper develops  the theory of
conformal invariants initiated in [\ref{R:F3}] for  Finsler
manifolds. More precisely we prove: The capacity of a compact set
and the capacity of the condenser of two closed sets are conformally
invariant. By mean of the notion of capacity, we construct and study
four conformal invariant functions $ \rho _{_M} $, $ \nu _{_M}$, $
\mu _{_M} $ and $\lambda _{_M}$ which have similarities with the
classical invariants on $S^n$, $ I\!\!R^n $ or $H^n$.
 Their properties and especially their continuity are
efficient tools for solving some problems of conformal geometry in
the large.
  \end{abstract}
\vspace{0.1cm}
  {\small{\textbf{Mathematics Subject Classification:} 30C70, 31B15, 53A30, 51B10.  }}
 \setcounter{page}{34}
 \thispagestyle{empty}
 \section*{Introduction.} The notion of
conformal capacity was introduced by Loewner [\ref{R:LO}] and has
been extensively developed for $ I\!\!R^n $ ( for instance
[\ref{R:G1}], [\ref{R:G2}], [\ref{R:MO}], [\ref{R:V2}]).
Particularly it was used by G.D Mostow to prove his famous theorem
on the rigidity of hyperbolic spaces [\ref{R:MO}]. J.Ferrand proved
that,  the capacity of  compact sets in Riemannian manifolds is
invariant under  conformal mappings and then she used this notion to
prove her famous theorem in Riemannian conformal  geometry
[\ref{R:F1}]. Here, inspiring her method,  we define an equivalent
notion of capacity in Finsler geometry and prove its invariance
property under conformal mappings.
\section{Preliminaries.}
\subsection{Finsler metric.}
 Let $M$ be a n-dimensional $ C^\infty$ manifold. For a point $x \in M$,
 denoted by $T_x M $
the tangent space of M at $x$. The tangent bundle of M is the union
of tangent spaces.
\[
TM:=\cup _{x \in M} T_x M.
\]
We will denote  the elements of TM by $(x, y)$ where $y\in T_xM$.
Let $TM_0 = TM \setminus \{ 0 \}.$ The natural projection $\pi: TM
\rightarrow  M$ is given by $\pi (x,y):= x$.
 Throughout this paper, we use {\it
Einstein summation convention} for the expressions with
indices.\\
 A {\it Finsler structure} on a manifold $M$ is a function $ F:TM_0
\rightarrow [0,\infty )$
with the following properties:\\
(i) $F$ is $C^\infty$ on $TM_0$.\\
(ii) $F$ is positively 1-homogeneous on the fibers of tangent
 bundle $TM$:
 \[
 \forall \lambda>0 \quad  F(x,\lambda y)=\lambda F(x,y).
 \]
 (iii) The Hessian of $F^{2}$ with elements
$ g_{ij}(x,y):=\frac{1}{2}[F^2(x,y)]_{y^iy^j} $ is positively
defined on $TM_0$. Then the pair  $(M,F)$ is called a {\it Finsler
manifold.} $F$ is Riemannian if $g_{ij} (x,y)$ are independent of $y
\neq 0$.
\subsection{Notations and definitions on conformal geometry
 of Finsler manifolds.}
 A diffeomorphism $f : (M,g)\rightarrow
(N,h)$ between $n$-dimensional
 Finsler manifolds $(M,g)$ and $(N,h)$
is called \textit{conformal}  if each $(f_{*})_p $ for $p\in M$ is
angle-preserving, and in this case two Finsler manifolds are called
\textit{conformally equivalent} or simply \textit{conformal}. If
$M=N$ then$f$ is called a \textit{conformal transformation
\textrm{or} conformal automorphism}. It can be easily checked that a
diffeomorphism is conformal if and only if \footnote{This result is
due to  Knebelman [\ref{R:KN}]. In fact the sufficient condition
implies that the function $\sigma(x,y)$ be independent of direction
$y$, or equivalently $\partial\sigma/\partial y^i = 0$.}, $f^{*}h =
e^{2\sigma} g$ for some function $\sigma: M\rightarrow \R$. The
diffeomorphism $f$ is called an \textit{isometry} if $ f^{*}h = g $.
 Now let's consider two Finsler manifolds $(M,g)$ and
$(\overline{M},\overline{g})$ with Finsler structures $F$ and $\bar
F$ and  with line elements $(x,y)$ and $(\bar{x},\bar{y})$
respectively. Throughout this paper we shall always assume that
coordinate systems on $(M,g)$ and $(\overline{M},\overline{g})$ have
been chosen so that  $ \overline{x^{i}}=x^{i}$ and $
\overline{y^{i}}=y^{i}$ holds, unless a contrary assumption is
explicitly made. Using this assumption we can show them by $(M,g)$
and $(M,\overline{g})$ or simply by $M $ and $\overline{M}$. Then
this two manifolds are conformal if $\overline{F}(x,y) =
e^{\sigma}F(x,y)$ or equivalently  $$\overline{g} =e^{2\sigma(x)}\
g\ .$$
 Locally we have $
\overline{g}_{ij}(x,y) = e^{2\sigma(x)}\ g_{ij}(x,y)$, and $
\overline{g}^{ij}(x,y) = e^{-2\sigma(x)}\ g^{ij}(x,y). $
\subsection{Some vector spaces and their properties.}
\subsubsection{Pull-back space $\pi ^* TM$.}
 Let $\pi : TM \longrightarrow M$
be the natural projection from $TM$ to $M$.\newline The {\it
pull-back tangent space  }$\pi ^* TM$ defined by
$$ \pi^* TM :=\{(x , y , v)|\,y \in T_x M_{0},\ v \in T_x M\}.$$
The {\it pull-back  cotangent space } $\pi^* T^*M$  defined by
$$\pi^* T^*M :=\{\pi ^*\theta |\,\theta \in T^*M\}.$$
Both $ \pi^* TM$ and  $\pi^* T^*M$ are n-dimensional vector spaces
over $T M_{0}$.\\

\subsubsection{Sphere bundle $SM$.}
 Let us denote by $S_xM$ the set consisting of all rays $[y]:=\{\lambda y
|\,\lambda > 0\},\,$ where $ y \in T_x M_{0}$ . Let
$$SM= \bigcup_{x \in M} S_xM.$$ $SM$ has a natural $(2n-1)$
dimensional manifold structure, called {\it Sphere bundle} over $M$.
We denote the elements of $SM$ by $(x , [y])$ where $y \in T_xM_{0}$
[\ref{R:BC}].\\
\begin{lem}{[\ref{R:BL}]}
The Sphere bundle of a differentiable manifold is orientable.
\end{lem}
\subsubsection{Pull-back space $ p^*TM$.}
Let $ p: SM \rightarrow  M$ denotes the natural projection from $SM$
to $M$. The {\it pull-back tangent space} $p^*TM$ is defined by
$$ p^*TM:= \{(x , [y] , v) |\, y \in T_xM_{0},\, v\in T_xM\}.$$
The {\it pull-back cotangent space} $p^*T^*M$ is defined by
$$ p^*T^*M:= \{p^*\theta |\,\theta\in T^*M\}.$$
Both $ p^*TM$ and $ p^*T^*M$ are n-dimensional vector spaces over $SM$.\\
Let we  define the function $\eta$ as follows $$ \eta : TM_0
\longrightarrow SM, $$
$$\eta (x , y) = ( x , [y]).$$
We use the following lemma for replacing the $ C^\infty$ functions
on $TM_0$ by those on $SM$.
 \bl{[\ref{R:NA}]} \label{L:21} Let $f
\in C^\infty (TM_0)$. Then there exist a function $g\in C^\infty
(SM)$ satisfying $\eta^*g = f$  if and only if $$f( x , y) = f (x,
\lambda y), \hspace {1cm} y\in T_xM_0 , \lambda>0,$$ where $\eta^*$
is the pull-back of $\eta$.
 \el

Let $f \in C^\infty (M)$, the vertical lift of $f$ is denoted by
$f^V \in C^\infty (TM_0)$ and defined by
$$f^V : TM \longrightarrow I\!\!R $$
$$ f^V ( x , y) := f \circ
\pi ( x , y) = f(x). $$
  $f^V$ is independent of  $y$ and from lemma \ref{L:21}  there is a
  function $g$
  on $C^\infty (SM)$ related to $f^V$ by means of $\eta^*g = f^V$.
  We denote $g$ in the sequel by $f^V$ for simplicity.

\subsection{ Nonlinear connection. }
\subsubsection{ On tangent bundle $TM$.}

  Consider $\pi_{*}:TTM \longrightarrow TM$ and let we put
$ker\pi_{*}^{v}=\{z\in TTM|\,\pi_{*}^{v}(z)=0 \},\,\, \forall v\in
TM,$  then the vertical vector bundle on $M$ is defined by $$VTM =
\bigcup_{_{v\in TM}}ker\pi_{*}^{v}.$$ A
 \textit{non-linear connection} or a \textit{horizontal distribution} on $TM$ is a
 complementary distribution $HTM$  for $VTM$ on $TTM$. The non-linear nomination arise from
 the fact  that $HTM$ is  spanned by the functions  which are completely
determined by the differentiable
  non-linear functions.
   These functions are called coefficients of the non-linear connection and will
  be noted in the sequel by $N^{j}_{i}$. It is clear that $HTM$ is a horizontal vector
  bundle.
 By definition we have the decomposition $TTM =VTM\oplus HTM$.\newline

 Using the induced coordinates $(x^{i},y^{i})$ on $TM$, where $x^{i}$ and $y^{i}$
 are called respectively {\it position} and {\it
direction} of a point on $TM$,  we have the local field of  frames
 $\{\frac{\partial}{\partial x_{i}},\frac{\partial}{\partial y_{i}}\}$ on $TTM$.
 Let $\{dx^i ,
 dy^i\}$
 be the dual  of $\{ \frac{\partial}{\partial x^i} ,  \frac{\partial}{\partial y^i} \}
  $. It
 is well known that we can choose a local field of frames
  $\{\frac{\delta}{\delta x^i},\frac{\partial}
 {\partial y_{i}}\}$ adapted to the above decomposition
  i.e. $\frac{\delta}{\delta x^i}\in {\cal X}(HTM)$ and
 $\frac{\partial}{\partial y_{i}}\in {\cal X}(VTM)$. They  are
  sections of  horizontal and vertical sub-bundle on
  $HTM$ and $VTM$,
  defined by
  $\frac{\delta}{\delta x^i}=\frac{\partial}{\partial x_{i}}-N^{j}_{i}\frac{\partial}
 {\partial y_{j}}$, where $N^{j}_{i}(x^{},y^{})$ are the coefficients of non linear
  connection.
 Clearly
$$ N^i_{\,\,j} = \gamma _{\,\,\,jk}^i y^k - C_{\,\,\,jk}^i\gamma
_{\,\,\,rs}^ky^ry^s,$$ where $ \gamma
_{\,\,\,jk}^i:=\frac{1}{2}g^{is}(\frac{\partial g_{sj}}{\partial
x^k} - \frac{\partial g_{jk}}{\partial x^s} + \frac{\partial
g_{ks}}{\partial x^j})$
 and
  $C_{ijk} = \frac{1}{2}\frac{\partial
g_{ij} }{\partial y^k}$. \\

\subsubsection{On sphere bundle  $SM$.}

Using the coefficients of non linear connection on $TM$ one can
define a non linear connection on $SM$ by using the objects which
are invariant under positive re-scaling $y \mapsto \lambda y$. Our
preference for being on $SM$ dictates us to work with
$$\frac{ N^i_{\,\,j}}{F} := \gamma _{\,\,\,jk}^i l^k -
C_{\,\,\,jk}^{\,i}\gamma _{\,\,\,rs}^kl^rl^s ,$$ where
$l^i=\frac{y^i}{F}.$ \\

 We prefer also to work with the local field of frames
$\{\frac{\delta}{\delta x^i} , F \frac{\partial}{\partial y^j}\}$
and $\{ dx^i , \frac{\delta y^j}{F} \}$ which are invariant under
the positive  re-scaling of $y$ and can be used as a local field of
frame for tangent bundle $ p^*TM$ and cotangent bundle $ p^*T^*M$
over $SM$ respectively.\newline

\subsection{A Riemannian metric on $SM$.}

 It turns out that the manifold $TM_0$ has a
natural Riemannian metric ( known in the literature as {\it Sasaki
metric} [\ref{R:BC}], [\ref{R:Mi}])
$$\widetilde{g} = g_{ij}(x, y)dx^i\otimes dx^j + g_{ij}(x, y)\frac{\delta
y^i}{F}\otimes\frac{\delta y^j}{F},$$ where $g_{ij}(x, y)$ are the
Hessian
 of Finsler structure $F^2$. They are functions on $TM_0$ and invariant under
positive re-scaling  of $y$, therefore they can be considered as
functions on $SM$. With respect to this metric, the {\it horizontal
subspace} spanned by $\frac{\delta}{\delta x^j}$ is orthogonal to
the {\it vertical subspace} spanned by $F \frac{\partial}{\partial
y^i}$.\newline
 The metric $\widetilde{g}$ is invariant
under the positive re-scaling of $y$ and can be considered  as a
Riemannian metric on $S(M)$.

\subsection{Hilbert form.}
Consider the pull-back vector bundle $p^*TM$ over $SM$. The
pull-back tangent bundle $p^*TM$ has a canonical section $l$ defined
by
$$ l_{(x , [y])}=(x , [y] , \frac{y}{F(x , y)}).$$
We use the local coordinate system $(x^i , y^i)$ for $SM$, where
$y^i$ being homogeneous coordinates up to a positive factor. Let
$\partial_i:=(x , [y] , \frac{\partial}{\partial x^i})$.
$\{\partial_i\}$ is a natural local field of frames for $p^*TM$. The
natural dual co-frame for $p^*T^*M$ is $\{dx^i\}$. The Finsler
structure  $F(x , y)$ induces a canonical 1-form on $SM$ defined  by
$$
\omega:= l_idx^i,
$$
where $ l_i=g_{ij}l^j.$\\
 $\omega$ is called {\it Hilbert form } of $F$. Using
 $g_{ij}= F F_{y^i y^j}+ F_{y^i}F_{y^j}$ and $\frac{\delta
F}{\delta x^i}=0 $, with straight forward calculation we get
\begin{equation}
\label{E:Bao}
 d\omega = -(g_{ij}- l_il_j)dx^i \wedge \frac{\delta
y^j}{F}.
\end{equation}

\subsection{Gradient vector field.}
  For a Riemannian manifold $ (S(M) , \widetilde{g})$,  the
gradient vector field of a function $f \in C^\infty (S(M))$  is
given by
$$ \widetilde{g}(\nabla f, \widetilde{X} )= df(\widetilde{X}),
\hspace{2cm} \forall \hspace{.2cm}\widetilde{X} \in {\cal X} (SM).$$
Using the local coordinate system $(x^i , [y^i])$ for $SM$, the
vector field $ \widetilde{X} \in {\cal X} (SM)$ is given by
$\widetilde{X} = X^i(x, y) \frac{\delta}{\delta x^i} + Y^i(x, y) F
\frac{\partial}{\partial y^j}$ where $X^i(x, y)$ and $Y^i(x, y)$ are
$C^\infty$ functions on $SM$. Using straight forward calculation we
get locally
$$ \nabla f =g^{ij} \frac{\delta f}{\delta x^i} \frac{\delta}{\delta
x^j} + F^2 g^{ij} \frac{\partial f}{\partial y^i} \frac{\partial}{
\partial y^j}.$$ The norm of  $ \nabla f $ with respect to the
Riemannian metric $\widetilde{g}$ is given by
\begin{equation}
\label{E:Norm}
 \mid\nabla f \mid ^{\,\, 2} =
\widetilde{g}( \nabla f , \nabla f) = g^{ij} \frac{\delta f}{\delta
x^i}\frac{\delta f}{\delta x^j}+ F^2 g^{ij} \frac{\partial
f}{\partial y^i}\frac{\partial f}{\partial y^j}.
\end{equation}
\section{Conformal invariant.}
In what follows $(M , g)$ denotes a  connected  Finsler manifold of
class $C^1$ and dimension $n \geq 2$. Let $(S(M) , \widetilde{g}) $
be its Riemannian Sphere
 bundle,
we set at first some definitions and notations.\\
Let's consider the Volume element $\eta (g) $ on $S(M)$ defined as
follows [\ref{R:AK}]
$$\eta (g) := \frac{(-1)^N}{(n-1)!}\,\, \omega \wedge (d \omega
)^{n-1},$$ where ${\small N=\frac{n(n-1)}{2}}$
and $\omega$ is a Hilbert form of $F$.\\

Let $H(M)=\mathcal{C}(M)\cap W^1_n (M)$ be the linear space of
continuous  real valued functions $u$ on $M$ admitting a generalized
$L^n$-integrable differential, satisfying
$$ I(u , M) = \int _{S(M)} \mid \nabla u^V\mid ^ n \eta (g) <\infty
,$$ where $u^V$ is the vertical lift of $u$.
\\

If $M$ is non-compact then $H_0(M)$ is the subspace of functions $u
\in H(M)$ such that its vertical lift  $u^V$ has a compact support
in $S(M)$.
\begin{defn} A function $u \in \mathcal{C}(M)$ will be called
monotone if for any relatively compact domain $D$ of $M$
$$\sup_{x \in \partial D} u(x)= \sup_{x \in D} u(x),
\,\,\,\,\,\,\,\,\,\,
 \inf_{x \in \partial D} u(x)= \inf_{x \in D}u(x).$$
We denote by $H^*(M)$ the set of monotone functions $u \in H(M)$.
\end{defn}
\begin{defn}
The capacity of a compact   subset  $C$ of a non-compact Finslerian
manifold $M$ is defined by
$$ Cap_{_M}(C) := \inf _u I(u , M), $$ where the infimum is taken
over the  functions  $u\in H_0(M)$ with $u=1$ on $C$ and $ 0\leq
u(x) \leq 1$ for all $x$, these functions being said to be
admissible for $C$.
\end{defn}
\begin{defn}Let $(C_0 , C_1)$ be a pair of closed
 sets in  Finslerian
manifold  $M$. The capacity of the condenser $\Gamma (C_0 , C_1, M)$
is defined by
 $$ Cap_{_M}(C_0 , C_1) = \inf _{u \in A(C_0 ,\, C_1)} I(u , M),$$
 where the infimum is taken over the set $A(C_0 , C_1)$ of all
 functions $u\in H(M)$ satisfying $u=0$ on $C_0$ and  $u=1$ on $C_1$
 and $ 0\leq u(x) \leq 1$ for all $x$, these functions are called
  admissible for condenser $\Gamma (C_0 , C_1, M)$.
If $A(C_0 , C_1) = \emptyset $ and particulary if $C_0 \bigcap C_1
\neq \emptyset $, we set
 $ Cap_{_M}(C_0 , C_1) = +\infty $.
\end{defn}
\begin{defn}
A relative continuum is a closed subset $C$ of $M$ such that $ C
\cup \{ \infty \}$ is connected in Alexandrov's compactification
$\overline{M}= M \cup \{ \infty \}$. For avoiding ambiguities the
connected closed sets of $M$ which are not reduced to one point will
be called continua.
\end{defn}
In what follows we want to associate the conformal invariant
functions determined entirely by conformal structure of manifold
$M$, at every double, triple and quaternary points of $M$.
\begin{defn} For all $(x_1 , x_2)$ in $M^2$ we set $$\mu_{_M}(x_1 , x_2)=
\inf_{C \in \alpha(x_1 , x_2)}Cap_{_M}(C),$$ where $\alpha(x_1 ,
x_2)$ is the set of all compact continua subsets of $M$, containing
$x_1$ and $ x_2$. And we set
$$\lambda_{_M}(x_1 , x_2)=\inf_{C_0 , C_1} Cap_{_M}(C_0 , C_1),$$
where $C_0$ and $C_1$ are relative continua resp. containing $ x_1 $
and $ x_2$.
\end{defn}
\begin{defn}Let $\triangle = \{ ( x , x , x)| \,\, x \in M\}$ be the
diagonal of $M^3$. For any $(x_1 , x_2 , x_3) \in M^3 \backslash
\triangle$ we set
$$\nu _{_M}( x_1 ,x_2 , x_3) = \inf_{C_0 , C_1} Cap_{_M}(C_0 , C_1),$$ where
$C_0$ is a relative continuum containing $x_3$ and $C_1$ a compact
continuum containing $x_1$ and $x_2$.
\end{defn}
\begin{defn}Let $\triangle$ be the set of all points $(x_1 ,x_2 , x_3
, x_4)$ of $M^4$ such that at least three coordinates of which  are
equal , and $\overline{I\!\!R}_+=I\!\!R_+ \cup  \{+ \infty\}$. We
define a function $\rho _{_M}: M^4 \backslash  \triangle
\longrightarrow \overline{I\!\!R}_+$ by setting $\rho _{_M} ( x_1
,x_2 , x_3 , x_4)= +\infty $ if $\{ x_1 , x_2\} \cap \{ x_3 ,
x_4\}\neq \emptyset $ and in all other cases
$$\rho _{_M} ( x_1 ,x_2 , x_3
, x_4)= \inf_{C_0 ,\, C_1} Cap_{_M}(C_0 , C_1),$$ where $C_0$ is a
compact continuum containing $ x_1, x_2$ and $C_1$ a compact
continuum containing $ x_3 , x_4$.
\end{defn}
\begin{defn} For any subset $S$ of $M$ and any $u \in
\mathcal{C}(M)$, we denote by $\omega (u , S)$ the oscillation of
$u$ on $S$.
\end{defn}

\section{Conformal properties of capacity.} Let $ f : M \longrightarrow M' $ be a
diffeomorphism between two manifolds and $h$ the mapping
$$h : S(M) \longrightarrow S(M'),$$
$$ h(x , [y])=(f(x) , [f_*(y)]),$$
 where $f_*$ is the differential of $f$ (the tangent map, [\ref{R:NA}]).
  Since $f_*$ is  linear, $h$ is well defined.
  Let $f$ be a conformal map between Finsler manifolds  $(M ,
g)$ and $ (M' , g')$, with the Finsler structures $F$ and $F'$
respectively. With respect to the function $\lambda$ on $M$ and
$\omega'$ be a Hilbert form related to the Finsler structure $ F'$.
In other word $ \omega'= g'_{ij}\frac{y'^j}{\sqrt{g'_{mn}y'^m
y'^n}}d x'^i$, we have
$$h^*\omega'= \sqrt{\lambda}\omega.$$
 from  (\ref{E:Bao}) we get $$ h^* d \omega'=
\sqrt{\lambda}d\omega.$$ So if $\eta (g) $ and $ \eta (g')$ denotes
the volume elements of $S(M)$ and $S(M')$ respectively, then we find
that
\begin{equation}
\label{E:OP}
 h^*( \eta (g'))= (\sqrt{\lambda})^n\eta (g).
\end{equation}
Therefore the mapping $h$ is  orientation preserving diffeomorphism
from $S(M)$ to $S(M')$. With above notions we have the following
lemma.
\begin{lem} If $u \in H_0 (M')$ then we have \newline

 I)  $\mid \nabla u^V \mid^n =( g'^{ij} \frac{\delta
u^V}{\delta x'^i}\frac{\delta u^V}{\delta
x'^j})^{\frac{n}{2}},$ \\

 II)  $(u o f)^V = u^V o h,$\\

 III) $h^* \frac{\delta u^V}{\delta x'^i}= \frac{\delta (u \circ f)^V}{\delta
 x^i}.$
\end{lem}
\begin{proof} The first assertion  follows from  (\ref{E:Norm}),  II)  and
III) can be easily verified  by direct calculations.
\end{proof}
From the above lemma we have
\begin{equation}
\label{E:hnorm} h^* \mid \nabla u^V \mid^n = (\sqrt{\lambda})^{-n}
\mid \nabla (u \circ f)^V \mid^n.
\end{equation}
Now we can prove the following theorem. It shows that, the capacity
of a compact  set and the capacity of the condenser of two closed
sets are conformally invariant, i.e. they only depend on the
conformal structure.
\begin{thm}
\label{T:1}Let $f$ be a conformal map between two  Finsler manifolds
$(M,g)$ and $(M',g')$. Then we have
$$Cap_{_M}(C)=Cap_{_{M'}}(f(C)),\,\,\,\,
 Cap_{_M}(C_0, C_1)=Cap_{_{M'}}(f(C_0), f(C_1)),$$
 for every compact subset $C$ and closed subsets $C_0$ and $C_1$ of
 $M$.
 \end{thm}
\begin{proof}
Let $f : (M , g) \longrightarrow (M' , g')$ be a conformal map.
First we prove \begin{equation} \label{L4:2} I(u , M')= I( u \circ f
, M), \end{equation}
 for every $u \in
H_0 (M')$. By definition
 $$ I( u , M') = \int _{S( M')} \mid
\nabla u^V \mid^n \eta (  g'),$$
Since $S(M)$ and $S(M')$ are two
orientable $n$-dimensional smooth manifolds with boundary and $h$ is
a smooth and orientation preserving diffeomorphism between them, we
have (see for example p. 245, [\ref{R:Lee}])
  $$ \int _{S( M')} \mid
\nabla u^V \mid^n \eta (  g') = \int _{S(M)} h^* ( \mid \nabla u^V
\mid^n \eta ( g')).$$ Using equation $(\ref{E:OP})$ and
$(\ref{E:hnorm})$ gives
$$\int _{S(M)} h^* ( \mid \nabla u^V
\mid^n \eta ( g')) = \int _{S(M)} \mid \nabla (u \circ f)^V\mid^n
\eta( g)= I( u \circ f , M).$$

 Let $C$ be a compact  set in $M$ by definition
$$ Cap_{_M}(C)= \inf_{v \in H_0M , \,\,v|_{_C}=1} I(v , M), \hspace{1cm}
   Cap_{_{M'}} (f(C))= \inf_{ u \in H_0M' ,\,\, u|_{_{f(C)}}=1} I(u ,
   M').$$
Putting $$A = \{ I(v , M) |v \in H_0M , v|_{_C}=1\},$$
    $$B = \{ I(u , M') | u \in H_0M' ,
    u|_{_{f(C)}}=1\},$$
    since $f^{-1}($ support $u)=$ support $(u\circ f)$ for all
$I(u , M') \in B,$ we have $(u \circ f) \in H_0(M).$ On the other
hand  $(u \circ f)|_{_C}=1$ and from relation (\ref{L4:2}), $I(u ,
M')= I( u \circ f , M)$.
 Hence $B \subseteq A $.\\
By the same argument we can prove $A \subseteq B $. Therefore $
Cap_{_M}(C)=Cap_{_{M'}}(f(C)).$ \\
 Let $C_0$ and $C_1$ be closed subsets of
 $M$. By  putting $$A = \{ I(v , M) |v \in H_0M ,
  v|_{_{C_0}}=0 , v|_{_{C_1}}=1\},$$
 $$B = \{ I(u , M') | u \in H_0M' ,
    u|_{_{f(C_0)}}=0 , u|_{_{f(C_1)}}=1 \},$$ with the same argument we
    can prove
     $Cap_{_M}(C_0,
C_1)=Cap_{_{M'}}(f(C_0), f(C1)).$
\end{proof}
By mean of the notion of capacity, we can study the properties of
four conformal invariant functions $ \rho _{_M} $, $ \nu _{_M}$, $
\mu _{_M} $ and $\lambda _{_M}$ which have similarities with the
classical invariants on $S^n$, $ I\!\!R^n $ or $H^n$
[\ref{R:F1}],[\ref{R:MO}].
 Their properties and especially their continuity are
efficient tools for solving some problems of conformal geometry.
\\In the following theorem we prove that the functions $\rho _{_M} ,\nu
_{_M} , \mu_{_M}$ and $\lambda_{_M}$ depend only on the conformal
structure of $M$ and therefore invariant under any conformal
mapping.
\begin{thm}

Let $f$ be a conformal mapping from the Finsler manifold $M$ to the
Finsler manifold $M'$ , we have for all $ x_1 , x_2 , x_3 , x_4
$ in $M$\\
$\rho _{_M}(x_1 , x_2 , x_3 , x_4 ) = \rho_{_{M'}}(f(x_1)  ,
f(x_2) , f(x_3), f(x_4)),$\\
$ \nu _{_M}(x_1 , x_2 , x_3 ) = \nu_{_{M'}}(f(x_1) , f(x_2) ,
f(x_3)),$\\
$\mu _{_M}(x_1 , x_2) = \mu_{_{M'}}(f(x_1)  , f(x_2)), $\\
$\lambda_{_M}(x_1 , x_2 ) = \lambda_{ _{M'}}(f(x_1) , f(x_2) ) .$
\end{thm}
\begin{proof}
The proof is a straight forward conclusion of theorem 1 and
definitions 5, 6 and 7.
\end{proof}
\newpage
{\small
{\it Author's address:}\\
B. Bidabad and S. Hedayatian\\
 Faculty of Mathematics, Amirkabir
University of Technology,\\ Tehran Polytechnic, 424, Hafez Ave.
15914, Tehran-Iran.\\
{\small\it E-mails:\textbf{ bidabad@aut.ac.ir} ${\,\,\,\,\,\,\,\,}$
\textbf{s$_{-}$hedayatian@aut.ac.ir}}}
\end{document}